\documentclass[12pt]{amsart}
\usepackage{amsmath}
\usepackage{amsfonts,amssymb,amsthm, txfonts, pxfonts}
\def\struckint{\mathop{%
\def\mathpalette##1##2{\mathchoice{##1\displaystyle##2}%
  {##1\textstyle##2}{##1\scriptstyle##2}{##1\scriptscriptstyle##2}}%
\mathpalette
{\vbox\bgroup\baselineskip0pt\lineskiplimit-1000pt\lineskip-1000pt
\halign\bgroup\hfill$}
{##$\hfill\cr{\intop}\cr\diagup\cr\egroup\egroup}%
}\limits}

\newtheorem{theorem}{Theorem}
\newtheorem{lemma}[theorem]{Lemma}
\newtheorem{corollary}[theorem]{Corollary}

\newtheorem{theorem-definition}[theorem]{Theorem-Definition}

\theoremstyle{remark}

\newcommand{\integers}{\mathbb{Z}}

\newcommand{\reals}{\mathbb{R}}

\DeclareMathOperator{\vol}{Vol}

\DeclareMathOperator{\mcg}{Mod}
\DeclareMathOperator{\stab}{Stab}

\begin{document}


\title{A simpler proof of Mirzakhani's simple curve asymptotics}

\author{Igor Rivin}

\address{Department of Mathematics, Temple University, Philadelphia}

\email{rivin@math.temple.edu}

\thanks{The author thanks Greg McShane for useful conversations, and
Benson Farb on helpful comments on a previous draft of this paper}

\date{\today}

\keywords{hyperbolic structures, mapping class group, ergodicity,
Mirzakhani, simple curves}

\subjclass{57M50, 32G15}

\begin{abstract}
Maryam Mirzakhani (in her doctoral dissertation) has proved the
author's conjecture that the number of simple closed curves of length
bounded by $L$ on a hyperbolic surface $S$ is \emph{asymptotic} to a
constant times $L^d,$ where $d$ is the dimension of the Teichm\"uller
space of $S.$ In this note we clarify and simplify Mirzakhani's argument.

\end{abstract}

\maketitle

\section*{Introduction} In \cite{mizphd}
M.~Mirzakhani shows that the number of closed geodesics of length
bounded above by $L$ in the mapping class group orbit of a fixed curve
$\gamma$ on a hyperbolic surface $X$ grows like a (positive) constant
times $L^{\dim T(X)},$ where $T(X)$ is the Teichm\"uller space of $X.$
This result, together with the obvious (but easy to miss) observation
that simple curves fall into a finite number of orbits, implies my
conjecture that the number $N_L(X)$  of simple closed geodesics of
length bounded above by $L$ has the same order of growth. Such a
result was obtain by G. McShane and the author in \cite{mc1,mc2}
in the special case where $X$ is homeomorphic to a punctured torus. In
the general case, I had previously obtained an order of growth result
in \cite{geomded}. The approach of \cite{geomded} was largely
combinatorial, and Mirzakhani has brought a number of different ideas
to bear on the question. These ideas are:

\begin{enumerate}
\item One should study the orbit of a fixed curve.
\item One can get a lot of H. Masur's result \cite{howieerg}: the
action of the mapping class group on the space of measured laminations
of a surface is \emph{ergodic}.
\item The growth rate is the product of two factors: one depends on
the hyperbolic metric on $X,$ while the other does not. One can then
``separate variables'' by integrating over moduli space.
\item The last integration uses Mirzakhani's result that the
Weil-Petersson volumes of moduli spaces of surfaces with boundary are
\emph{polynomials} in the length of the boundary.
\end{enumerate}

In this note, we combine these ideas with the results of
\cite{geomded}, the coarea formula, and some simple observations to
simplify Mirzakhani's argument.

\section{Multicurves}
\label{multicurves}
The first part of the argument concerns multicurves on a hyperbolic
surface $S.$ A \emph{multicurve} is a collection of pairwise disjoint simple
closed curves on $S,$ together with the assignment of an integer
weight to each curve in the collection (the weight is usually thought
of as the number of times the curve winds around itself). 

There is a way to parametrize the set of all multicurves (discovered
by Max Dehn): we decompose $S$ (in some fixed way) into pairs of
pants, then record, for each pants boundary curve, how many times the
multicurve intersects it and how many times it winds around
it. Winding can be clockwise or counterclockwise, while intersection
numbers are \emph{geometric} intersection numbers, and so are
positive (and also even). It is an easy result of Dehn's that there is
a bijective corespondence between the collections of such integer
coordinates and multicurves; the reader can consult Dehn's original
paper \cite{dehn} or the more recent \cite{harerpenner} for the proofs
of the assertions above. 

Dehn's coordinates thus represent multicurves
as integer lattice points in an improper cone in $\reals^{\dim
\mathcal{ML}(S)}.$\footnote{For a closed surface of genus $g,$ the
dimension of the moduli space (and of the measured lamination space)
is $6g - 6,$ for a surface with $c$ punctures and $p$ perforations,
the dimension is $6g - 6 + 2 c + 2 p.$ Since none of the arguments
below depend on the signature of the surface, we will always denote
the dimension by $\dim \mathcal{ML}(S)$ or $\dim \mathbb{M}(S),$ (the
two are equal) as
appropriate.}
This cone can be thought of as the \emph{measured lamination} space $\mathcal{ML}$ --
that identification is the content of Theorem 3.1.1 in
\cite{harerpenner}.
Explicitely 
\[
\mathcal{ML}(S) = \reals^{(\dim \mathcal{ML})/2} \times \reals_+^{(\dim
\mathcal{ML})/2},
\]
where the first half of the coordinates correspond to twists, and the
second to the intersection numbers; multicurves correspond to points
in $\mathcal{ML}(S)$ where the first half of the coordinates are
integral, and the second half are even.

 On this cone we have a measure, induced from
Lebesgue measure in the ambient $\reals^n.$ As is well known to
physicists far and wide, there is another way to define the Lebesgue
measure of an open set $\Omega$ in $\reals^n:$
\[
\lambda(\Omega) = \lim_{t\rightarrow \infty} \dfrac{\left| t \Omega
\cap \integers^n\right|}{t^n}.
\]
This is just Riemann integration in $\reals^n.$ For future reference,
let us also define
\[
\lambda_t(\Omega)  = \dfrac{\left| t \Omega \cap
\integers^n\right|}{t^n},
\]
so that 
\[
\lambda(\Omega) = \lim_{t\rightarrow \infty} \lambda_t(\Omega).
\]
It follows from these definitions that the Lebesgue measure is
invariant under the Mapping Class Group of $S,$ (henceforth denoted by
$\mcg(S),$ since $\mcg(S)$ acts bijectively on integer lattice points
and respects scaling.

\section{Length functions}
Dehn's coordinates are essentially topological, but a hyperbolic
structure on a surface $S$ defines a \emph{length function} on
multicurves, where the length of a multicurve is simply the weighted
sum of the lengths of the geodesic representatives of the connected
components of the multicurve in question. Indeed, one way to think of
the length is as the inf of the lengths of topological (multi)curves
isotopic to the given multicurve. It is then clear that the length
function is linear on rays (this follows from the uniqueness part of
Dehn's theorem alluded to above) and otherwise is convex (since
$L(a+b) \leq L(a) + L(b)$). The length function can thus be extended
by linearity to all points with rational coordinates, and by
continuity (which follows from the triangle inequality above) to all
real points in $\mathcal{ML}$ (This construction is identical to and
slightly more general than the construction in \cite{mc1,mc2},
as is the immediate sequel). From the convexity of the length function
$L$ it follows that the set
\[
B_L(1) = \{x \left| L(x) < 1 \right.\}
\]
is a convex set in $\mathcal{ML}.$ Further, by linearity and the
definition of the measure $\lambda$ in Section \ref{multicurves}, it
follows that
\begin{theorem}
\label{multiasymp}
The number of multicurves of length bounded above by $L$ is asymptotic
to 
\[
\lambda(B_L(1)) L^{\dim \mathcal{ML}}.
\]
\end{theorem}
Theorem \ref{multiasymp} would be useless, unless we knew that
\[
0 < \lambda(B_L(1)) < \infty.
\]
Luckily, that is well-known, and follows from the fact that the $L^1$
norm of the Dehn coordinate vector is quasi-the-same as the hyperbolic
length of the corresponding curve (the distortion depends on the
hyperbolic structure, and diverges as the hyperbolic structure goes to
the boundary of Teichm\"uller space). Essentially that result is
explained in \cite{geomded}, though there is no doubt tha the result
was known to experts for at least twenty years beforehand.

An additional observation is that $\lambda(B_L(1))$ varies
analytically over moduli space. This follows, eg, from the methods of
\cite{kerckan}.

\section{Orbit density}
\label{orbitdensity}
Consider a set $X \subseteq \mathcal{ML}$ and a multicurve $x \in
\mathcal{ML}.$ Let $O(x)$ denote the orbit of $x$ under the mapping
class group, and define
\[
\mu_t(X) = t^{-\dim \mathcal{ML}} \left|O(x) \cap t X \right|.
\]
Each $\mu_t$ defines a $\mcg(S)$-invariant measure, and it is clear
that the family $\{\mu_t\}$ is bounded, thus weakly compact. Since
$\mu_t(X) \leq \lambda_t(X)$ for all $t,$ it follows that for any
subsequence $\sigma,$ the subsequence limit $\mu_\sigma$ is absolutely
continuous with respect 
to $\lambda.$ Since the action of $\mcg(S)$ on $\mathcal{ML}$ is
\emph{ergodic,} it follows that such a subsequence limit $\mu_\sigma$
is a constant multiple of $\lambda.$ That is, for any finite
measurable $X,$ 
\[
\mu_\sigma(X) = c_\sigma(x) \lambda(X).
\]
In particular, we can take $X = B_L(1)$ for a length function $L$
coming from a hyperbolic structure on $S.$ 
The results of
\cite{geomded} give the following
\begin{theorem}
\label{rivconst}
There exists a \emph{curve} $x,$ and a constant $c$ such that 
$c_\sigma(x) > c$ for \emph{any} subsequence $\sigma.$
\end{theorem}
\begin{proof} I show in \cite{geomded} that the number of simple
curves of length bounded by $L$ has order of growth $L^{\dim
\mathcal{ML}}.$ Since there is only a finite number of $\mcg(S)$
orbits of simple closed curves, at least one of the orbits grows at
that speed, and this gives the desired $x.$
\end{proof}

\section{Curves and moduli}
Take a surface $S$ and a \emph{curve}\footnote{There is no
need for $\gamma$ to be a curve -- all of the arguments work
\emph{mutatis mutandis} for a disjoint collection of essential
curves. The argument does not work directly for general multicurves,
though the result for such is easily deduced.} $\gamma.$
Take the cover $\mathbb{M}^\gamma$ of the moduli space $\mathbb{M}(S)$
corresponding to the subgroup $\stab(\gamma) \subset \mcg(S).$
If we cut $S$ along $\gamma,$ we obtain a (possibly disconnected)
surface $S^\prime$ with two boundary components of the same length,
and it is clear that this induces a map $\pi_\gamma :
\mathbb{M}^\gamma \rightarrow \mathbb{M}_\gamma S^\prime,$ where the
last subscript indicates that two of the boundary components
correspond to the same curve. The map is a fibration, where the fiber
corresponds to twisting along $\gamma.$ If we take a pair of pants
decomposition of $S$ which includes $\gamma$ as one of the pants
curves, $\pi$ acts very simply on the Weil-Petersson symplectic form:
the term $d \ell(\gamma) \wedge d \tau(\gamma)$ is simply dropped. The fiber is a
circle of length $L(\gamma).$ This shows (by the coarea formula
(Corollary \ref{coarea}) that the volume of that piece of
$\mathbb{M}^\gamma$ where the length of the appropriate translate of
$\gamma$ equals $l$ is simply $l$ times the volume of the moduli space
of $S^\prime$ where the two curves corresponding to $\gamma$ have
length $l.$ This is so because Wolpert's formula (see \cite{wolpertsymp}) implies that the
Jacobian of $\pi_\gamma$ equals $1.$

We remind the reader that Wolpert's formula for the Weil-Petersson
K\"ahler form reads:
\[
WP = \sum_{\gamma} d \ell(\gamma) \wedge d \tau(\gamma),
\]
where the sum is taken over all curves in a pants decomposition of
$S,$ $\ell(\gamma)$ is the length of $\gamma$ and $\tau(\gamma)$ is
the twist along $\gamma.$

Let us now denote the number of curves in the $\mcg(S)$ orbit of
$\gamma$ of length not exceeding $L$ (corresponding to a hyperbolic
structure $\mathcal{H}$) by $n_{\mathcal{H}}(L),$ and let us ask what
the \emph{average} value of $n_{\mathcal{H}}(L)$ is over the moduli
space of $S.$ It turns out to be easier to not normalize by the
volume of $\mathbb{M}(S).$ In that case, we have the following obvious
relationship:
\begin{equation}
\label{obviousintegral}
\int_{\mathbb{M}(S)} n_{\mathcal{H}(L)} = \vol\left\{(x, \rho)  \in
\mathbb{M}^\gamma \left|\quad L(\rho) < L \right.\right\},
\end{equation}
where the notation $(x, \rho)$ means that we are in the sheet
corresponding to the image $\rho$ of $\gamma.$ The volume in the right
hand side of the Eq. \eqref{obviousintegral} is easy to evaluate with
the help of the discussion at the beginning of the section, the coarea
formula, and Mirzakhani's results on the Weil-Petersson volumes of
moduli spaces of bordered surfaces: the answer is a polynomial of
degree $\dim \mathbb{M}(S) = \dim \mathcal{ML}.$ This gives us the
following observation:
\begin{lemma}
\label{isconstant}
There exists a constant $C,$ such that 
\[
L^{-\dim \mathbb{M}(S)} \int_{\mathbb{M}(S)} n_{\mathcal{H}(L)} < C.
\]
\end{lemma}
Now we can show that the constant $c_\sigma(x)$ in Section
\ref{orbitdensity} does not depend on the subsequence $\sigma.$
Indeed, with the limit being taken over an arbitrary subsequence
$\sigma:$
\[
\lim \int_{\mathbb{M}(s)} n_{\mathcal{H}(L)} L^{-\dim(\mathbb{M}(S))} =
\int_{\mathbb{M}(S)} \lim n_{\mathcal{H}(L)} L^{-\dim(\mathbb{M}(S))},
\] by the Dominated Convergence Theorem and Lemma
\ref{isconstant}. The left hand side does not depend on the
subsequence by the discussion following
Eq. \eqref{obviousintegral}. The right hand side equals
\[
\int_{\mathbb{M}(S)} c_\sigma(\gamma) \lambda(B_L(1)),
\]
and so also does not depend on the choice of the subsequence $\sigma.$

Note that if we pick $\gamma$ to be the curve whose existence is shown
in Theorem \ref{rivconst}, it follows that:
\begin{theorem}
\label{finiteness}
The integral of $\lambda(B_L(1))$ over $\mathbb{M}(S)$ is finite.
\end{theorem}
It then follows immediately that $c_\sigma(\gamma)$ is itself
independendent of $\sigma.$

\section{The coarea formula}
In this setting we need a very simple version of the coarea formula of
Federer, but we shall state a more general version. Our source for
this is Ralph Howard's exposition \cite{howardcoarea}. First, $Jf$ is
the \emph{Jacobian} of $f,$ defined (under the assumption that $m \geq
n$) as:
\[
Jf(x) = \begin{cases}
0, & \text{if $x$ is a critical point of $f,$}\\
\sqrt{\det(f_*(x) f_*^t(x))}, & \text{otherwise}.
\end{cases}
\]

\begin{theorem}[The Coarea Formula]
\label{gencoarea}
 Let $f: M^m \rightarrow N^n$ be a
smooth map between Riemannian manifolds, with $m \geq n.$ Then, for
almost every $y \in N^n,$ the fiber $f^{-1}(y)$ is either empty or a
submanifold of $M^m$ of dimension $m - n.$ For each regular value $y$
of $f,$ let $dA$ be the $m - n$ dimensional surface area measure on
$f^{-1}(y).$ Then, for any measurable function $h$ on $M^m,$
\[
\int_{N^n}\int_{f^{-1}(y)} h d A d Y = \int_{M^m} h(x) Jf(x) d x,
\]
where $dy$ is the Riemannian volume measure on $N^n,$ and $dx$ is the
Riemannian volume measure on $M^m.$
\end{theorem}
\begin{corollary}
\label{coarea}
Let $\mathcal{H}^{m-n}(f^{-1}(y))$ denote the $m-n$ dimensional area
of $f^{-1}(y).$ Then
\[
\int_{N^n} \mathcal{H}^{m-n}(f^{-1}(y)) d y = \int_{M^m} Jf(x) d x.
\]
\end{corollary}
\begin{proof}
Set $h \equiv 1$ in Theorem \ref{gencoarea}.
\end{proof}
\bibliographystyle{plain}
\bibliography{amie,rivin}

\begin{thebibliography}{10}

\bibitem{dehn}
Max Dehn.
\newblock Lecture notes from breslau.
\newblock Technical report, 1922.
\newblock Archive of the University of Texas at Austin.

\bibitem{howardcoarea}
Ralph Howard.
\newblock The kinematic formula in {R}iemannian homogeneous spaces.
\newblock {\em Mem. Amer. Math. Soc.}, 106(509):vi+69, 1993.

\bibitem{kerckan}
Steven~P. Kerckhoff.
\newblock Earthquakes are analytic.
\newblock {\em Comment. Math. Helv.}, 60(1):17--30, 1985.

\bibitem{howieerg}
Howard Masur.
\newblock Ergodic actions of the mapping class group.
\newblock {\em Proc. Amer. Math. Soc.}, 94(3):455--459, 1985.

\bibitem{mc2}
Greg McShane and Igor Rivin.
\newblock Geometry of geodesics and a norm on homology.
\newblock {\em International Mathematics Research Notices}, (2):61--69,
  February 1995.

\bibitem{mc1}
Greg McShane and Igor Rivin.
\newblock Simple curves on hyperbolic tori.
\newblock {\em Comptes Rendus Acad. Sci. Paris, S\'er. I. Math},
  320(12):1523--1528, June 1995.

\bibitem{mizphd}
Maryam Mirzakhani.
\newblock Growth of the number of simple closed geodesics on hyperbolic
  surfaces.
\newblock submitted, 2004.

\bibitem{harerpenner}
R.~C. Penner and J.~L. Harer.
\newblock {\em Combinatorics of train tracks}, volume 125 of {\em Annals of
  Mathematics Studies}.
\newblock Princeton University Press, Princeton, NJ, 1992.

\bibitem{geomded}
Igor Rivin.
\newblock Simple curves on surfaces.
\newblock {\em Geometriae Dedicata}, 87(1/3):345--360, August 2001.

\bibitem{wolpertsymp}
Scott Wolpert.
\newblock On the symplectic geometry of deformations of a hyperbolic surface.
\newblock {\em Ann. of Math. (2)}, 117(2):207--234, 1983.

\end{thebibliography}
\end{document}